\documentclass[a4paper,11pt,reqno]{amsart}
\usepackage{graphicx}
\usepackage{color}
\definecolor{rblue}{RGB}{39,64,139}
\usepackage[pdftex]{hyperref}
\hypersetup{
	naturalnames=false,%
	draft=false,
	allcolors=rblue,
	allbordercolors={.15 .25 .55},
	colorlinks=true,
	}
\usepackage{amssymb}
\usepackage[initials]{amsrefs}
\usepackage[cmtip]{xypic}

\title[Borel-Weil-Bott Theorem via Equivariant McKean-Singer Formula]{Borel-Weil-Bott Theorem via\\ Equivariant McKean-Singer Formula}
\author{Seunghun Hong}
\address{Mathematisches Institut, Busenstra{\ss}e 3--5, D-37073 G\"ottingen, Germany}
\email{shong@uni-goettingen.de}
\urladdr{diracoprerat.org}

\keywords{Borel-Weil-Bott theorem, cubic Dirac operator, equivariant index theorem, equivariant McKean-Singer formula}
\subjclass[2010]{Primary %
19K56
, 58J35
, 43A77
; Secondary %
43A85
, 22E45
, 58A14
}

\thanks{This research was partially supported under NSF grant DMS-1101382.}

\theoremstyle{plain}
\newtheorem{thm}{Theorem}

\theoremstyle{definition}

\theoremstyle{remark}
\newtheorem*{rmk}{Remark}

\newcommand{\KD}{\mathcal{D}}
\newcommand{\C}{\mathbb{C}}

\newcommand{\Z}{\mathbb{Z}}

\newcommand{\grg}{\mathfrak{g}}

\newcommand{\grp}{\mathfrak{p}}

\newcommand{\spin}{\mathfrak{spin}}

\newcommand{\Id}{\mathbf{1}}
\DeclareMathOperator{\Spin}{Spin}
\DeclareMathOperator{\SO}{SO}
\DeclareMathOperator{\Str}{Str}
\DeclareMathOperator{\tr}{tr}

\DeclareMathOperator{\Ind}{Ind}

\DeclareMathOperator{\Cl}{Cl}
\DeclareMathOperator{\Ad}{Ad}
\DeclareMathOperator{\diag}{diag}

\newcommand{\spinor}{\mathbb{S}}
\DeclareMathOperator{\Hom}{Hom}

\begin{document}

\begin{abstract}
After reviewing how the Borel-Weil-Bott theorem can be interpreted as an index theorem, we present a proof using Kostant's cubic Dirac operator and the equivariant McKean-Singer formula. 
\end{abstract}
\maketitle

The Borel-Weil-Bott theorem ``completes'' the representation theory of compact Lie groups by providing a cohomological method to construct each and every irreducible representation of a compact Lie group. The main constituents of this construction are the harmonic forms, that is, the homogeneous solutions to the Hodge-Dolbeault operator. Thus the Borel-Weil-Bott method is essentially an index theorem. This viewpoint underlies the work of Bott~\cite{bott} and is manifest in Slebarski's theorem \cite{slebarski1}*{Thm.~1, p.~296}. In this article, we give a brief review on the interpretation of the Borel-Weil-Bott theorem as an index theorem, and prove it using Kostant's cubic Dirac operator and the equivariant McKean-Singer formula.

\section{Preliminaries}\label{sec:prelim}
Throughout the article we write $X=G/T$, where $G$ is a compact connected Lie group and $T$ is a maximal torus of $G$. The Lie algebra of $G$, that is, the tangent space of $G$ at the identity $e$, shall be denoted by the lowercase black letter $\mathfrak{g}$. We endow $\mathfrak{g}$ with an inner product $\langle \cdot , \cdot \rangle$ by taking the negative of the Killing form. This inner product is invariant under the adjoint action of $G$ on $\mathfrak{g}$.

\subsection{Some Basic Facts and Notations Surrounding the Representation Theory of Compact Lie Groups}

Consider the conjugation action of $G$ on itself. The elements among $G$ that preserve $T$ constitute the normalizer $N_G(T)$ of $T$ relative to $G$. The quotient 
\[ W=N_G(T)/T\]
is known as the \emph{Weyl group} of $G$. 

The action of $W$ on $T$ induces a $W$-action on functions on $T$; in particular, the irreducible characters of $T$. The irreducible characters of $T$ constitute the unitary dual $\hat T$ of $T$. So the orbit space $\hat T/W$ makes sense. A consequence of the \emph{Weyl character formula} is that there is a one-to-one correspondence between $\hat T/W$ and the unitary dual $\hat G$ of $G$.

Let $\theta\colon T\to \C^\times$ be an irreducible character of $T$. Its derivative $\theta_*\colon \mathfrak{t}\to \C$ is a Lie algebra representation. The function $-i\theta_*$ is a linear functional on the real vector space $\mathfrak{t}$. The image  of the map $\hat T\to \mathfrak{t}^*$, $\theta\mapsto -i\theta_*$, forms a lattice $\Lambda_T$ in $\mathfrak{t}^*$. We have:
\[ \Lambda_T = \{\, \lambda\in\mathfrak{t}^*\,|\,\lambda(H)\in2\pi\Z\ \text{for all} \ H\in\mathfrak{t}\cap\exp^{-1}\{e\}\,\},\]
where $\exp$ denotes the exponential map. The constituents of $\Lambda_T$ are known as the \emph{analytically integral weights} of $G$. 

Let $K$ denote the fundamental Weyl chamber of our choice for the $W$-action on $\mathfrak{t}^*$. Then the one-to-one correspondence between $\hat G$ and $\hat T/W$ implies a one-to-one correspondence between $\hat G$ and $\Lambda_T\cap K$. More precisely:
\[\begin{array}{ccc} \hat G  &\leftrightarrow & \Lambda_T\cap K,\\
 {[V]} &\mapsto & \text{highest weight of $V$}.
 \end{array}
\]
For each $\lambda\in \Lambda_T\cap K$, we denote by $V_\lambda$ an irreducible $G$-representation space according to the above correspondence.

Let $w\cdot \lambda$ denote  the action of $w\in W$ on $\lambda\in\Lambda_T$. The  \emph{shifted action} of $w$ on $\lambda$ is:
\[ w\odot \lambda := w\cdot(\lambda+\rho)-\rho,\]
where $\rho$ is the Weyl vector, that is, half the sum of the positive roots. (The positive roots are determined by the Weyl chamber $K$.) Then  $\Lambda_T\cap K$ consists of all elements of $\Lambda_T$ that has free $W$-orbit with respect to the shifted $W$-action. In summary, the Weyl character formula implies the following one-to-one correspondence:
\begin{equation}
\begin{array}{ccc} \hat G  &\leftrightarrow & \{\text{free shifted $W$-orbits in $\Lambda_T$}\},\\
 {[V_\lambda]} &\mapsto & W\odot \lambda.
 \end{array}
\label{eq:weylcharoto}
\end{equation}

\subsection{Borel-Weil-Bott Theorem}
Let $\mu\in \Lambda_T\cap K$, and let $ U_\mu$ denote the complex vector space $\C$ on which $T$ acts by the irreducible character of weight $\mu$. We denote by $\ell(\mu)$ the word length of $\mu$ relative to the fundamental Weyl chamber $K$; it satisfies:
\[ \ell(\mu)=\#\{\,\alpha\in\Phi_+ \mid  \left<\mu,\alpha\right> <0\,\}. \]
Here $\Phi_+$ denotes the set of positive roots of $G$, and $\langle \cdot , \cdot \rangle$ is the inner product on $\mathfrak{g}^*$ induced by that on $\mathfrak{g}$.

Let $G\times_T U_\mu$ denote the space of equivalence classes in $G\times U_\mu$ with respect to the relation $(g,z)\sim(gx^{-1},x\cdot z)$ for $x\in T$. This is a complex line bundle over $X$. It is diffeomorphic to
\[ \mathcal{L}_\mu:=G^{\C}\times_B U_\mu,\]
which is a complex line bundle over $G^{\C}/B$, where $G^{\C}$ is the complexification of $G$, and $B$ is a Borel subgroup of $G^{\C}$. 

The statement of the Borel-Weil-Bott theorem depends on the selected Borel subgroup $B$. Our convention is as follows. Let $\mathfrak{g}_{\C}$ be the complexification of $\mathfrak{g}$, and take the root space decomposition:
\[ \mathfrak{g}_{\C}= \mathfrak{n}_-\oplus\mathfrak{h}\oplus\mathfrak{n}_+,\]
where $\mathfrak{n}_\pm$ denotes the positive and negative root spaces, and $\mathfrak{h}$ is the complexification of $\mathfrak{t}$. We set
\[ \mathfrak{b}:= \mathfrak{h}\oplus\mathfrak{n}_+.\]
Then $B$ is the connected subgroup in $G^{\C}$ with Lie algebra $\mathfrak{b}$.

Let $\mathcal{O}(\mathcal{L}_\mu)$ be the sheaf of germs of holomorphic sections of $\mathcal{L}_\mu$. The celebrated Borel-Weil-Bott theorem states that the sheaf cohomology $H^*(X;\mathcal{O}(\mathcal{L}_\mu))$ is nontrivial only if $\mu$ has free shifted $W$-orbit, and if that is the case then
\[ H^q(X;\mathcal{O}(\mathcal{L}_\mu)) \cong \begin{cases}
	V_{W\odot\mu},& q=\ell(\mu);\\
	0, &\text{otherwise}. \end{cases}
	\]
Here $V_{W\odot \mu}$ denotes the irreducible representation space of $G$ corresponding to the shifted orbit of $\mu$ according to the correspondence~\eqref{eq:weylcharoto}.

Now consider the twisted Hodge-Dolbeault complex
\[  \mathcal{A}^p_\mu := \Omega^{0,p}(X)\otimes U_\mu,\]
whose differential is given by the Dolbeault operator
\[ \bar\partial := d^{0,1}\otimes\Id.\]
Owing to the Dolbeault theorem (see \cite{wells}*{Thm.~3.20, p.~63}), the complex $(\mathcal{A}^\bullet_\mu,\bar\partial)$ computes the sheaf cohomology of $\mathcal{O}(\mathcal{L}_\mu)$:
\[ H^*(X;\mathcal{O}(\mathcal{L}_\mu))\cong H^*\{(\mathcal{A}^\bullet_\mu,\bar\partial)\}.\]

Meanwhile, by the Hodge theorem, $H^*\{(\mathcal{A}^\bullet_\mu,\bar\partial)\}$ is isomorphic to the kernel of the Dirac operator 
\begin{equation}
 D:=(\bar\partial +\bar\partial^\dag)/\sqrt2,
 \label{eq:dolbeultdirac}
\end{equation}
where $\bar\partial^\dag$ is the formal adjoint of $\bar\partial$. Since $D$ is $G$-equivariant, the kernel of $D$ is a $G$-representation space; we denote the corresponding virtual representation as $[\ker D]$. Then Borel-Weil-Bott theorem is equivalent to saying: 
\begin{equation}
  [\ker D] =\begin{cases}
  	[V_{W\odot\mu}], & \text{if $W\odot\mu$ is a free orbit;}\\
	0,&\text{otherwise};
	\end{cases}
\label{eq:bwbequiv}
\end{equation}
and $\ker D$ is homogeneous in degree equal to $\ell(\mu)$. This form of the Borel-Weil-Bott theorem first appeared in Slebarski~\cite{slebarski1}.

\subsection{Borel-Weil-Bott Theorem as an Equivariant Index Theorem}
The complex $\mathcal{A}^\bullet_\mu$ is naturally bi-graded by the even and odd forms. Let $D_+$ and $D_-$ denote the restrictions of $D$ onto the even and odd subspaces. We have
\[ [\ker D]= [\ker D_+]+[\ker D_-].\]
Since $V_{W\odot \mu}$ is irreducible (when $W\odot \mu$ is free), Equation~\eqref{eq:bwbequiv} can be refined as follows. If $W\odot\mu$ is free then:
\[
  [\ker D] =[\ker D_+] \quad\text{or}\quad [\ker D]=[\ker D_-].
\]
If $W\odot \mu$ is not free then:
\[
 [\ker D]=[\ker D_+]=[\ker D_-]=0.
 \]

Now the \emph{equivariant index} of $D$ is by definition the virtual representation
\[ [\Ind D]:=[\ker D_+]-[\ker D_-].\]
Owing to what we have just seen above, we have:
\[ [\Ind D] = [\ker D_+] = [\ker D] \quad\text{or}\quad [\Ind D] = -[\ker D_-] = -[\ker D],\]
provided that $W\odot\mu$ is free; otherwise we have:
\[ [\Ind D] = [\ker D_\pm]=[\ker D]=0.\]
Thus, the Borel-Weil-Bott theorem implies that $[\Ind D]$ is nontrivial if and only if $\mu$ has free shifted $W$-orbit, and if that is the case then $[\Ind D]$ is equal to $[V_{W\odot\mu}]$ up to sign. In fact, in our proof of the Borel-Weil-Bott theorem, we shall show that: 
\begin{equation}
  [\Ind D] = \begin{cases}
  (-1)^{\ell(\mu)}[V_{W\odot \mu}],&\text{if $W\odot\mu$ is a free orbit};\\
  0,&\text{otherwise}.
  \end{cases}
\label{eq:bwbeqvind}
\end{equation}
Equation~\eqref{eq:bwbequiv} then follows with the aid of Equation~\eqref{eq:flgmfldind2}. This index theorem is a refinement of Bott's result \cite{bott}*{Thm.~III, p.~170}. It was first shown by Landweber~\cite{landweber} (for the general case of compact homogeneous space $G/H$ where $H$ is a closed subgroup of maximal rank in $G$).

\subsection{Equivariant McKean-Singer Formula and Kostant's Cubic Dirac Opeartor}
In obtaining Equation~\eqref{eq:bwbeqvind}, Landweber uses Bott's equation: $[\Ind D] = i_*([E]-[F])$, where $E$ and $F$ denotes $T$-spaces, and $i_*$ is the induction map $R(T)\to \hat{R}(G)$, $[E]\mapsto [\Gamma^2(G\times_TE)]$. Our method is to use, in place of Bott's equation, the equivariant McKean-Singer formula. What we have in mind more precisely is this: A virtual representation can be identified with its image under the character map 
\[ \chi\colon R(G)\to C(G), \]
which maps an irreducible element $[V]$ to its character $\chi_V$. For the value of $\chi_V$ at  $g\in G$, we write   
\[ [V]_g := \chi_V(g). \]
The equivariant McKean-Singer formula then states that $[\Ind D]_g$ is equal to the super trace of the operator $ge^{t\KD^2}$ where $t$ is a positive real number (see Berline, Getlzer, and Vergne \cite{bgv}*{Prop.~6.3, p.~185}):
\begin{equation}
 [\Ind D]_g = \Str(ge^{tD^2})=\tr(ge^{tD_-D_+}) -\tr(ge^{tD_+D_-}). \label{eq:equivmcksing} 
\end{equation}
(Although the right-hand side seems at first to be an infinite linear combination of irreducible characters, it is actually a finite combination due to the symmetry between the eigenvalues of $D_+$ and $D_-$.)
 
Instead of directly working with the complex $\mathcal{A}^\bullet_\mu$ in calculating the super trace, we shall use the isomorphism:
\begin{equation}
 \mathcal{A}^\bullet_\mu \cong (C^\infty(G)\otimes\wedge^\bullet(\mathfrak{n}_+)\otimes U_\mu)^T.
 \label{eq:dolcplxiso}
\end{equation}
Here the action of $T$ is as follows: On $ U_\mu$ it is by the irreducible character of weight $\mu$; on $C^\infty(G)$ it is the one induced by right-translations; and on $\wedge^\bullet(\mathfrak{n}_+)$ it is that induced by the adjoint action. The $T$-action on $\wedge^\bullet (\mathfrak{n}_+)$ is related to the spinors constructed out of the orthogonal complement $\mathfrak{p}$ of $\mathfrak{t}$ in $\mathfrak{g}$ by:
\begin{equation}
 [\wedge^\bullet(\mathfrak{n}_+) ]  = [\spinor^*\otimes U_{\rho}] \in  R(T).
 \label{eq:nilspinid}
\end{equation}
Here $\spinor^*$ is dual of the spinor space $\spinor$ associated to the Clifford algebra $\Cl(\mathfrak{p})$ generated by $\mathfrak{p}$ (see Kostant~\cite{kostant2000}*{Prop.~3.6, p.~76}). The action of $T$ on $\spinor$ is provided by taking the homomorphism $T\to \SO(\mathfrak{p})$, coming from the adjoint representation of $G$, and lifting it to  $T\to \Spin(\mathfrak{p})$:
\[ \xymatrix@C=4em{ &\Spin(\mathfrak{p})\ar[d]\\
	T\ar[ru]^{\widetilde{\Ad}} \ar[r]_-{\Ad} &\SO(\mathfrak{p})}\]
This lift always exists \cite{freedflag}*{Cor.1.12, p.~91}. In short, we have:
\begin{equation}
 \mathcal{A}_\mu^\bullet\cong (C^\infty(G)\otimes \mathbb{S}^*\otimes  U_{\mu+\rho})^T.
\label{eq:dolcplxiso2}
\end{equation}

Finally, because the equivariant index of a Dirac operator depends only on its symbol (this is easy to check directly, but there is a general theorem by Bott~\cite{bott}*{Thm.~I, p.~169}), we may use, in place of the Dirac operator~\eqref{eq:dolbeultdirac}, Kostant's cubic Dirac operator:
\begin{equation} \KD := \sum_{i=1}^{\dim\mathfrak{p}} Y_i\otimes Y_i +1\otimes \frac{1}{3} \sum_{i=1}^{\dim\mathfrak{p}}Y_i\gamma(Y_i)\in \mathcal{U}(\mathfrak{g}) \otimes \Cl(\mathfrak{p}). \label{eq:kostantdirachmsp}
\end{equation}
Here $\mathcal{U}(\mathfrak{g})$ is the universal enveloping algebra of $\mathfrak{g}$; $\{Y_i\}_{i=1}^{\dim\mathfrak{p}}$ is any orthonormal basis for $\mathfrak{p}$; and $\gamma$ is the map $\mathfrak{g}\to \spin(\mathfrak{p})$ defined by:
\begin{equation}
\gamma(Z):=-\frac{1}{2}\sum_{i,j=1}^{\dim\grp}\langle Z,[Y_i,Y_j]_\grg\rangle Y_iY_j.\label{eq:decmopsospin2}
\end{equation}
The action of the algebra $\mathcal{U}(\mathfrak{g})\otimes\Cl(\mathfrak{p})$  on the right-hand side of~\eqref{eq:dolcplxiso2} is trivial on $ U_{\mu+\rho}$; the action on $\spinor^*$ comes from the canonical action of $\Cl(\mathfrak{p})$; and the action on $C^\infty(G)$ is solely from $\mathcal{U}(\mathfrak{g})$, which arises from identifying $Z\in\mathfrak{g}$ with the left-invariant vector field it generates on $G$.

The advantage of using the cubic Dirac operator lies in the simple form of its square:
\begin{equation}
 \KD^2  = -\Omega_\mathfrak{g}+\diag\Omega_\mathfrak{t} + \|\rho\|^2.
 \label{eq:kostdircsq}
\end{equation}
Here $\Omega_\mathfrak{g}$ denotes the Casimir element in $\mathcal{U}(\mathfrak{g})$, and $\diag$ denotes the algebra homomorphism $\mathcal{U}(\mathfrak{t})\to \mathcal{U}(\mathfrak{g})\otimes \Cl(\mathfrak{p})$ induced by the map $\mathfrak{t}\to \mathcal{U}(\mathfrak{g})\otimes \Cl(\mathfrak{p})$, $X\mapsto X\otimes 1+1\otimes\gamma(X)$. 

To see the  effectiveness of Equation~\eqref{eq:kostdircsq}, decompose the right-hand side of~\eqref{eq:dolcplxiso2} using the Peter-Weyl theorem:
\begin{equation}
 (C^\infty(G)\otimes\spinor^*\otimes U_{\mu+\rho})^T \cong\bigoplus_{[V_\lambda]\in\hat G} V_\lambda\otimes (V_\lambda^*\otimes \spinor^*\otimes  U_{\mu+\rho})^T.
 \label{eq:ptwytwstcplx}
\end{equation}
(This is not entirely correct; the isomorphism is true upon taking the norm closures on both sides.) This isomorphism is obtained by identifying ${\left|v\right\rangle}{\otimes}{\left\langle w\right|}\in V_\lambda\otimes V_\lambda^*$ with the function $G\to \C$, $g\mapsto \left\langle g\cdot w\,|\, v\right\rangle$. The action of $\mathcal{U}(\mathfrak{g})$ on the right-hand side is the one induced by the Lie algebra representation on each $V_\lambda^*$. Now take a summand 
\[ H_\lambda:=V_\lambda\otimes (V_\lambda^*\otimes \spinor^*\otimes  U_{\mu+\rho})^T.\]
Owing to Schur's lemma, the action of $\Omega_\mathfrak{g}$ on $V_\lambda^*$ is constant with the value  $-\left\|\lambda+\rho\right\|^2+\|\rho\|^2$. For similar reasons, the action of $\diag\Omega_\mathfrak{t}$ on $V_\lambda^*\otimes\spinor^*$ is again constant with the value  $-\|\mu+\rho\|^2$. In summary, the restriction of $\KD^2$ on the summand $H_\lambda$ is simply the constant operator
\begin{equation}
 \KD^2_\lambda := \|\lambda+\rho\|^2-\|\mu+\rho\|^2.
 \label{eq:kosconst}
\end{equation}
So the super trace of the operator $ge^{-t\KD^2}$ restricted to $H_\lambda$ is:
\begin{equation}
 \Str(ge^{-t\KD^2_\lambda}) = [V_\lambda]_g\left< [V_\lambda\otimes \spinor_+]-[V_\lambda\otimes\spinor_-],[ U_{\mu+\rho}]\right>_Te^{-t(\|\lambda+\rho\|^2-\|\mu+\rho\|^2)},
\label{eq:strkdtwstcplx}
\end{equation}
where $\langle \cdot , \cdot \rangle_T$ denotes the nondegenerate paring on $R(T)$ defined by:
\[ \left<E,F\right>_T = \dim \Hom_T(E,F). \]

\section{A Proof of the Borel-Weil-Bott Theorem via the Equivariant McKean-Singer Formula}

We now derive the Borel-Weil-Bott theorem using the equivariant McKean-Singer formula. As we have explained in Section~\ref{sec:prelim}, the Borel-Weil-Bott theorem is equivalent to the following:

\begin{thm}\label{thm:bwb} Let $\KD$ be Kostant's cubic Dirac operator acting on the smooth sections of the twisted vector bundle $G\times_T(\spinor\otimes U_{\mu+\rho})$ over $X=G/T$. The equivariant index of $\KD$ satisfies: 
\begin{equation}
 [\Ind\KD] = \begin{cases}
	(-1)^{\ell(\mu)}[V_{W\odot \mu}],&\text{if $W\odot\mu$ is free};\\
	0,&\text{otherwise}.
	\end{cases}
\label{eq:flgmfldind}
\end{equation}
Moreover,
\begin{equation}
[\Ind\KD]  = 
\begin{cases}
\phantom{-}[\ker \KD_+], & \text{if $\ell(\mu)$ is even};\\
-[\ker \KD_-], & \text{if $\ell(\mu)$ is odd}.\\
\end{cases}
\label{eq:flgmfldind2} \end{equation}
In each case, contributions to $\ker \KD_\pm $ comes from sections whose degree is equal to $\ell(\mu)$.
\end{thm}
\begin{proof}
We begin by invoking the equivariant McKean-Singer formula:
\[
[\Ind \KD]_g = \sum_{[V_\lambda]\in\hat G} \Str(ge^{t\KD^2_\lambda}). \]
By Equation~\eqref{eq:strkdtwstcplx}, we have:
\[ [\Ind \KD]_g= \sum_{[V_\lambda]\in\hat G} [V_\lambda]_g\left< [V_\lambda\otimes \spinor_+]-[V_\lambda\otimes\spinor_-],[ U_{\mu+\rho}]\right>_Te^{-t(\|\lambda+\rho\|^2-\|\mu+\rho\|^2)}.
\]
But the left-hand side is independent of the parameter $t$; hence, the only contribution in the sum occurs from the terms with the exponential factor equal to $1$, that is, when $\|\lambda+\rho\|=\|\mu+\rho\|$. Thus we have:
\[
[\Ind \KD] = \sum_{\substack{[V_\lambda]\in\hat G,\\ \|\lambda+\rho\|=\|\mu+\rho\| }}[V_\lambda]\left< [V_\lambda\otimes \spinor_+]-[V_\lambda\otimes\spinor_-],[ U_{\mu+\rho}]\right>_T.
\]
According to the multiplicity result of Kostant~\cite{kostant}*{Thm.~4.17, p.~486}, we have:
\begin{equation}
\Bigl(\bigl< [V_\lambda\otimes \spinor_+],[ U_{\mu+\rho}]\bigr>_T,\bigl<[V_\lambda\otimes\spinor_-],[ U_{\mu+\rho}]\bigr>_T\Bigr)  = 
\begin{cases}
(1,0),& \text{if $\ell(\mu)$ is even},\\
(0,1),& \text{if $\ell(\mu)$ is odd},
\end{cases}
\label{eq:kostmult}
\end{equation} 
provided that $\mu\in W\odot \lambda$; this last condition can be satisfied by some $[V_\lambda]\in \hat G$ if and only if $W\odot \mu$ is free. As a consequence we have Equation~\eqref{eq:flgmfldind}. We also find from Equation~\eqref{eq:kostmult} that a nontrivial contribution to $[\Ind \KD]$ comes solely from the even or the odd domain according to the parity of $\ell(\mu)$; hence Equation~\eqref{eq:flgmfldind2} holds. The same multiplicity result of Kostant also implies that such contribution to $[\Ind\KD]$ comes from elements whose degree is $\ell(\mu)$. This completes the proof.
\end{proof}

\begin{rmk}
Theorem~\ref{thm:bwb} can be modified so that it holds for more general cases where $T$ may be any closed subgroup $H$ of $G$ that is of maximal rank; the only change necessary is that we replace $U_{\mu+\rho}$ with $U_{\mu+\rho'}$, where 
\[ \rho':=\frac12\sum_{\alpha\in \Phi_+\smallsetminus\Phi_+(\mathfrak{h})}\alpha.\]
Here $\Phi_+(\mathfrak{h})$ denotes the set of positive roots of the Lie algebra $\mathfrak{h}$ of $H$ (the roots are calculated with respect to a common maximal toral subalgebra of $\mathfrak{g}$ and $\mathfrak{h}$). This change is necessary because Equation~\eqref{eq:nilspinid} now takes the form:
\[ [\wedge^\bullet(\mathfrak{n}_+) ]  = [\spinor^*\otimes U_{\rho'}].\]
The formula for $\mathcal{D}^2$ also changes to:
\[ \KD^2  = -\Omega_\mathfrak{g}+\diag\Omega_\mathfrak{h} + \|\rho\|^2 - \|\rho_\mathfrak{h}\|^2, \]
where $\rho_\mathfrak{h}$ is the Weyl vector of $\mathfrak{h}$. But Equation~\eqref{eq:kosconst} remains unmodified; so the argument we gave for $G/T$ can be repeated word-for-word for $G/H$, and we have the full results of Landweber~\cite{landweber}*{Thm.~3, p.~471} and Slebarski~\cite{slebarski2}*{Thm.~2, p.~509}.  
\end{rmk}

\end{document}